# Enhanced Energy Management System with Corrective Transmission Switching Strategy — Part I: Methodology

Xingpeng Li, *Student Member, IEEE* and Kory W. Hedman, *Member, IEEE*

*Abstract*— The flexibility in transmission networks is not fully utilized in existing energy management systems (EMSs). Corrective transmission switching (CTS) is proposed in this two-part paper to enable EMS to take advantage of the flexibility in transmission systems in a practical way. This paper proposes two EMS procedures: 1) Procedure-A connects real-time security-constrained economic dispatch (RT SCED) with real-time contingency analysis (RTCA), which is consistent with industrial practice; 2) Procedure-B, an enhanced version of Procedure-A, includes CTS in EMS with the proposed concept of branch pseudo limit used in RT SCED. Part-I of this paper presents the methodology while Part-II includes detailed results analysis. It is demonstrated that Procedure-A can effectively eliminate the base case overloads and the potential post-contingency overloads identified by RTCA and Procedure-B can achieve significant congestion cost reduction with consideration of CTS in RT SCED. Numerical simulations also illustrate that integrating CTS into RT SCED would improve social welfare.

## NOMENCLATURE

*Sets*
| | |
|---|---|
| $C$ | Contingencies. |
| $D$ | Loads. |
| $D(n)$ | Loads at bus $n$. |
| $DN$ | Negative loads. |
| $DV$ | Virtual loads. |
| $G$ | Generators. |
| $G(n)$ | Generators at bus $n$. |
| $GD$ | Generators that are available for dispatch. |
| $I$ | Interfaces. |
| $IKM(0)$ | Interface lines under monitor for base case. |
| $IKM(c)$ | Interface lines under monitor for contingency $c$. |
| $IM(0)$ | Interfaces under monitor for base case. |
| $IM(c)$ | Interfaces under monitor for contingency $c$. |
| $K$ | Branches. |
| $K(n-)$ | Branches of which bus $n$ is the from-bus. |
| $K(n+)$ | Branches of which bus $n$ is the to-bus. |
| $KI(i)$ | Branches that form interface $i$. |
| $KM(0)$ | Branches under monitor in base case. |
| $KM(c)$ | Branches under monitor for contingency $c$. |
| $N$ | Buses. |

*Parameters*
| | |
|---|---|
| $BS_{g,i}$ | Breadth of segment $i$ of unit $g$. |
| $C_{g,i}$ | Cost for segment $i$ of unit $g$. |
| $CSR_g$ | Spinning reserve price of unit $g$. |
| $C_{sh}$ | A fixed penalty factor for load shedding. |
| $EC$ | Operating energy cost in SCED. |
| $LimitA_k$ | Normal limit in MW of branch $k$ in SCED. |
| $LimitC_k$ | General emergency limit in MW of branch $k$ for all contingencies in SCED. |
| $LimitC_{kc}$ | Customized emergency limit in MW of branch $k$ under contingency $c$ in SCED. |
| $Limit_i$ | Total flow limit of interface $i$ in base case. |
| $Limit_{i,c}$ | Total flow limit of interface $i$ for contingency $c$. |
| $LL_k$ | Loading level of branch $k$ in base case. |
| $LL_{kc}$ | Loading level of branch $k$ under contingency $c$. |
| $LODF_{k,c}$ | Branch $c$ outage distribution factor on branch $k$. |
| $MRR_g$ | Energy ramp rate of unit $g$. |
| $NS_g$ | Number of cost segments for unit $g$. |
| $OTDF_{n,k,c}$ | Outage transfer distribution factor from bus $n$ to branch $k$ when branch $c$ is forced to be out of service. |
| $Pct$ | Branch monitoring tolerance for base case constraints. |
| $PctC$ | Branch monitoring tolerance for contingency case constraints. |
| $P_d$ | Forecasting load at the end of a SCED period. |
| $P_{d0}$ | Initial load at the beginning of a SCED period. |
| $P_{g0}$ | Initial output of unit $g$. |
| $P_{g,max}$ | Maximum output of unit $g$. |
| $P_{g,min}$ | Minimum output of unit $g$. |
| $P_{k0}$ | Initial MW flow on branch $k$ in SCED. |
| $P_{k,c,0}$ | Initial MW flow on branch $k$ under contingency $c$. |
| $PI_{n0}$ | Initial injection of bus $n$. |
| $PTDF_{n,k}$ | Power transfer distribution factor from bus $n$ to branch $k$. |
| $RateA_k$ | Normal limit in MVA of branch $k$. |
| $RateC_k$ | Emergency limit in MVA of branch $k$. |
| $RC$ | Operating reserve cost in SCED. |
| $SRR_g$ | Spinning ramp rate of unit $g$. |
| $T_{ED}$ | Look-ahead time of SCED. |
| $T_{SR}$ | Response time for spinning reserve requirement. |
| $x_k$ | Reactance of branch $k$. |
| $\alpha_k$ | Phase shifting angle of branch $k$; 0 for non-phase-shifter branch. |
| $nb$ | Number of buses. |
| $n_{(k-)}$ | From-bus of branch $k$. |
| $n_{(k+)}$ | To-bus of branch $k$. |

*Variables*
| | |
|---|---|
| $p_{d,sh}$ | Shedded active power of load $d$. |

This work was supported by 1) the Department of Energy Advanced Research Projects Agency - Energy, under the Green Electricity Network Integration program and under the Network Optimized Distributed Energy Systems program, and 2) the National Science Foundation award (1449080).

Xingpeng Li and Kory W. Hedman are with the School of Electrical, Computer, and Energy Engineering, Arizona State University, Tempe, AZ, 85287, USA (e-mail: Xingpeng.Li@asu.edu; kwh@myuw.net).

| | |
|---|---|
| $p_g$ | Output of unit $g$. |
| $p_{g,i}$ | Output on segment $i$ of unit $g$. |
| $p_k$ | Power flow on branch $k$. |
| $p_{k,c}$ | Power flow on branch $k$ under contingency $c$. |
| $sr_g$ | Spinning reserve that unit $g$ provides. |
| $\delta_n$ | Phase angle of bus $n$. |
| $\delta_{n,c}$ | Phase angle of bus $n$ under contingency $c$. |
| $PI_n$ | Power injection of bus $n$. |
| $AvgLMP$ | Average LMP over all buses in the system. |
| $AvgLMP_{cg}$ | Average congestion component of LMP. |
| $CCR_{CTS}$ | Congestion cost reduction due to CTS. |
| $CngstCost$ | Total congestion cost. |
| $CngstRvn$ | Total congestion revenue. |
| $GenCost$ | Total generator cost. |
| $GenRent$ | Total generator rent. |
| $GenRvn$ | Total generator revenue. |
| $LdPaymt$ | Total load payment. |
| $LMP_n$ | LMP at bus $n$. |
| $LMP_s$ | System-wide LMP. |
| $LMP_{cg,n}$ | Congestion component of LMP at bus $n$. |

## I. INTRODUCTION

As electric power cannot be economically stored on a large scale, it must be produced, transferred, and consumed at the same time. This creates serious challenges to maintain reliable real-time operations of power systems. Thus, an energy management system (EMS), a computer-aided tool, is used to help system operators monitor, control, and optimize real-time operations of electric power systems.

Some key functions of EMS include real-time contingency analysis (RTCA) and real-time security-constrained economic dispatch (RT SCED). RTCA evaluates the impact of a potential contingency on system security while RT SCED aims to provide a least-cost dispatch solution that meets the operation and reliability requirements. The transmission network is modeled as static assets in most (if not all) existing EMSs; in other words, the flexibility in transmission networks is not modeled and utilized in existing real-time operational tools. However, previous research has demonstrated that treating transmission network as a flexible asset can benefit the system in various aspects. In addition, operators can temporarily reconfigure the network topology in practice [1]-[4].

Corrective transmission switching (CTS) is considered as an efficient and practical strategy to take advantage of the transmission network flexibility for relieving the post-contingency violations identified by RTCA. A benchmark is provided in [5] to show that network topology control can improve system reliability by reducing post-contingency voltage violations and branch flow violations. Heuristic CTS algorithms are proposed in [6]-[7] to provide fast quality switching solutions for handling contingency-induced violations; over 1.5 million contingencies simulated on three large-scale real power systems demonstrate the proposed CTS heuristics. An enhanced data mining method is proposed in [8] to find the beneficial CTS solutions with much less time than other effective approaches. A comparison between dynamic CTS methods and a static look-up table (PJM's switching solutions) shows that dynamic CTS methods have much better performances [9]. Overloads and voltage violations caused by contingencies can be relieved by line and bus-bar switching actions [10]. Though [5]-[10] show that system reliability can be improved by incorporating CTS into RTCA, it is unclear how other EMS applications such as RT SCED that use RTCA results as input information can take advantage of the reliability benefits provided by CTS.

Optimal transmission switching (OTS) determines the best transmission topology and generation dispatch points simultaneously, which can reduce the cost significantly. For instance, it is found that a saving of 25% is achieved on the IEEE 118-bus system [11]; though cost saving drops to 15% when $N$-1 reliability is considered, it is still substantial [12]. OTS is typically formulated as a mixed integer linear programming (MILP) problem that is time consuming to solve [13]-[15]. Thus, heuristics are developed to reduce the computing time [14]-[15]. A tractable algorithm is proposed to solve the OTS problem and numerical simulations demonstrate the proposed tractable algorithm achieved 93% of the potential cost savings [14]. Two heuristics are proposed to reduce the solution time by iteratively solving a sequence of relaxed problems [15]. It is shown in [16] that topology control can significantly reduce congestion management costs. Network reconfiguration is considered as an effective tool for congestion management [17]. Corrective transmission topology control is shown to have the capability of improving reserve deliverability [18]. Integrating topology control into unit commitment problem can improve reserve deliverability, relieve congestion, and reduce cost [19]-[22]. Though [11]-[22] demonstrates various benefits by incorporating transmission network flexibility into scheduling and dispatching applications, they are based on the simplified DC model and their AC feasibility is not verified. To resolve this issue, some research that combine AC optimal power flow (AC OPF) model with OTS has been performed [23]-[25]. However, they are either taking too much computing time (~3,000 seconds for a 118-bus system) [23] or do not meet reliability requirements such as $N$-1 criterion [24]-[25]. It is worth noting that in reality, AC OPF has not been implemented yet after its first formulation was proposed over 50 years ago in 1962 [26], not mentioning the even more complex problem combining AC OPF with OTS. In addition, none of [11]-[25] investigated how to connect transmission switching based optimal scheduling and dispatching with RTCA.

Though prior efforts in the literature have demonstrated a variety of benefits by treating transmission network as a flexible asset in short-term operational applications [5]-[25], the information exchange and interfaces between different applications or modules have been ignored. To bridge this gap, this work proposes a seamless connection strategy between RTCA and RT SCED that can practically utilize the flexibility in the transmission network. One main challenge for implementing transmission switching for real-time operations is computational complexity. The status of switchable branches is typically represented by binary variables but that comes at the cost of a long solution time, which is impractical for real-time operations. To resolve this concern, pseudo limit is proposed in this work to avoid the use of binary variables in RT SCED.

In this two-part paper, corrective transmission switching is proposed to enable RTCA based SCED to take advantage of the flexibility in the transmission systems in a practical way. First, a traditional EMS procedure that mimics existing industrial practice is proposed in this paper to connect AC based RTCA and DC based RT SCED; this procedure is referred to as Procedure-A. A novel EMS procedure that considers the reliability benefits provided by CTS is presented to enhance the proposed Procedure-A; this enhanced procedure is referred to as Procedure-B in this paper. In addition, five different SCED models are proposed and evaluated, and the most effective model is selected for further examining the performance of Procedure-A and Procedure-B.

Fig. 1 illustrates the real-time operations of electric power systems. The key functions of EMS include state estimation, RTCA, and RT SCED. Measurements from remote terminal units (RTUs) are collected by state estimation through the communication network. Then, state estimation will process the measurements and estimate the system status, which will provide a base case for the subsequent EMS applications. As shown in Fig. 1(a), for the proposed Procedure-A, after system status is determined by state estimation, RTCA will scan the system and send system vulnerability information in the form of critical contingencies and the associated violations to RT SCED; then, RT SCED will redispatch online generators to eliminate those potential violations. As shown in Fig. 1(b), for the proposed Procedure-B, an additional CTS module is added in the EMS between RTCA and RT SCED to practically utilize the flexibility in transmission network and achieve both reliability benefits and economic benefits.

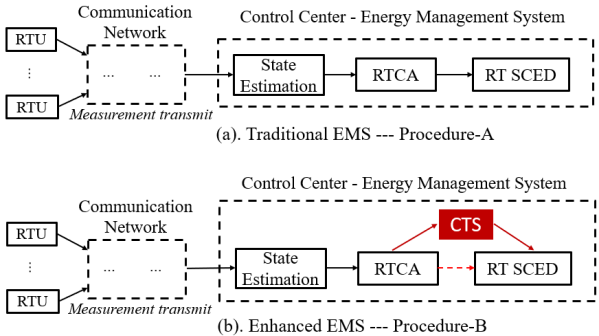

Fig. 1. Illustration of the power grid energy management system.

The main contributions of this paper are summarized below:
(i) The proposed Procedure-A practically connects AC based RTCA and DC based SCED, which respects the industrial practice. It is shown that Procedure-A can efficiently relieve the actual violations and the potential post-contingency overloads identified by RTCA.
(ii) The proposed concept, branch pseudo limit, enables SCED to avoid the use of binary variables when integrating CTS into SCED. This makes the SCED problem with consideration of CTS a linear programming (LP) problem, which substantially reduces the computational complexity.
(iii) With the proposed pseudo limit, the proposed procedure-B can be implemented with minimal change in existing EMS tools of electric power systems.
(iv) The proposed Procedure-B enables the enhanced SCED (E-SCED) to take advantage of the flexibility in transmission networks practically by using the proposed branch pseudo limit of network constraint rather than the actual limit for a traditional SCED.
(v) This paper investigated how much economic benefits can be translated from the reliability benefits provided by CTS. It is shown that with Procedure-B, CTS can achieve both substantial violation reduction and significant cost savings.
(vi) It is observed the CTS solution that reduces branch overloads for a contingency in the pre-SCED stage can also provide benefits for the same contingency in the post-SCED stage.
(vii) The energy market results associated with Procedure-B and E-SCED demonstrate that integrating CTS in EMS can reduce the system overall LMP and lower the total cost and thus improve the social welfare. Note that the proposed two EMS procedures can also be applied to the regulated power systems.

The rest of this paper is organized as follows. Section II presents a literature review on RTCA, RT SCED, and CTS. Section III describes the two proposed EMS procedures. Detailed RT SCED models are introduced in Section IV. Finally, Section V concludes the paper.

## II. LITERATURE REVIEW

### A. Real-Time Contingency Analysis

The North American Electric Reliability Corporation (NERC) requires system operators to maintain delivery of electric power after a single contingency or the loss of a single element ($N$-1) [27]. To ensure system $N$-1 reliability, RTCA is performed on a regular basis in real-time. It identifies the critical contingencies that may cause network violations, which helps operators to be prepared in advance and react to critical outages by using pre-determined strategies.

Though various classes of generation reserve requirements are enforced in power system scheduling and dispatching, $N$-1 reliability is not guaranteed. Thus, to better comply with NERC $N$-1 reliability standard, independent system operators (ISOs) conduct RTCA successively every few minutes while the actual implementation of RTCA varies for different ISOs.

A physical power system can be represented by the node-breaker-branch model which contains more details than the bus-branch model. However, to reduce the problem dimension and improve the algorithm robustness, a node-breaker-branch model will be converted into a simplified bus-branch model via graph theory techniques before any EMS applications including RTCA are performed.

Generators' outputs are assumed to remain the same for branch contingency while a participation factor based generation redispatch is performed following a generator contingency [5]-[6]. Due to the low probability of occurrence, multiple contingencies are not simulated in this work, which is in line with the industrial practice.

### B. Real-Time Security-Constrained Economic Dispatch

RT SCED is an optimization process that aims to provide the least cost solution while meeting all the physical, operation, and reliability constraints. It executes repeatedly every few minutes with updated real-time information. It typically runs automatically and upon operators' demand.

Though AC power flow model based SCED is accurate, it is not implemented in practice since AC SCED is a non-linear and non-convex problem and is extremely difficult to solve for large-scale real power systems in a timely manner. In addition, convergence is also of big concern. Hence, the simplified DC model based SCED is preferred and implemented at ISOs. However, the actual RT SCED models implemented at different ISOs are not the same. For instance, the RT SCED tool used at California ISO covers multiple look-ahead intervals while PJM performs RT SCED on a single interval only.

There are two types of SCED, corrective SCED and preventive SCED [28]. Corrective SCED allows the units' outputs under contingency to deviate from their base-case generation schedule, which would reduce the cost. The base-case solution determined by preventive SCED is supposed to withstand contingency situations without any adjustments. Thus, preventive SCED can provide a more secure and reliable solution. Due to reliability concerns, preventive SCED is typically implemented in practice rather than corrective SCED. Thus, to be consistent with industry practices, the SCED implemented in this work is preventive SCED.

It is worth mentioning that in the SCED models implemented by ISOs, all constraints are relaxed with slack variables. Thus, the SCED will never terminate due to infeasibility even in the worst scenario; instead, it can inform operators of the sources causing constraint violations if any.

Non-convexity such as generator prohibited operating zones may be involved in the SCED application; they must be addressed to ensure SCED solutions do not violate any physical restrictions. Binary variables can be used to directly model the generator prohibited operating zone in SCED; however, this will significantly increase the computational complexity of SCED by converting it from a LP problem into a MILP problem. Thus, heuristics are often implemented to avoid non-convexity in SCED. For instance, LP based SCED can be solved iteratively until the solution meets all physical constraints: if current iteration leads to a SCED solution that the dispatched power of a unit is within a prohibited operating zone; then, a constraint will be added to SCED to prevent that generator's output power from being dispatched within the same prohibited operating zone again in the next iteration. Since the proposed procedures does not change existing SCED model, the same heuristics can be applied to the proposed procedures. Note that our work is focused on SCED rather than the heuristics that deal with special non-convexities. Thus, non-convexity is beyond the scope of this work and is not considered in this work.

Theoretically, it is possible that the optimizer for SCED does not solve within an acceptable time. However, practically, this rarely happens because SCED does not involve any binary variables and it is a linear programming problem. Thus, SCED is typically solved very fast. The various SCED models presented in this paper are all LP problems, which is consistent with industrial practice.

### C. Corrective Transmission Switching

Corrective transmission switching switches a transmission element out of service shortly after a contingency occurs as a corrective control to reduce post-contingency violations. Thus, it can be used to relieve network congestion in RT SCED.

Prior work in the literature has demonstrated various benefits by treating transmission network as a flexible asset, including reliability enhancement [5]-[10], cost saving [11], [20], [29], losses reduction [14]-[15], congestion management [17]-[19], [30], and management of renewable uncertainty [31]-[33]. However, the flexibility in transmission networks is not fully utilized in contemporary power system real-time operations where the transmission network is modeled as a fixed topology. In this work, corrective transmission switching is proposed to take advantage of the flexibility in transmission network in the EMS functions including RT SCED.

### III. EMS PROCEDURES

In the proposed Procedure-A, the network constraints enforced in SCED are determined by RTCA; model conversion is performed to connect DC model based SCED with the traditional AC RTCA. This is consistent with the industrial practice. Extended upon Procedure-A, Procedure-B is proposed to utilize the flexibility in transmission networks by connecting SCED with CTS-based RTCA.

Procedure-A uses the actual limits for network constraints determined by RTCA while higher pseudo limits are used in Procedure-B, which is a main difference between Procedure-A and Procedure-B. The pseudo limits are determined by CTS-based RTCA. With the use of higher pseudo limits, the extra reliability provided by CTS can be considered in SCED and be translated into economic benefits.

Theoretically, CTS can be directly modeled in SCED with binary variables indicating the status of switchable element, which will convert SCED from a LP problem into a MILP problem. This will cause serious computational burden and substantially increase the solution time. Therefore, directly modeling CTS in SCED is impractical and Procedure-B with pseudo limits is preferred. In addition, the proposed Procedure-B requires no change to existing SCED tools and the corresponding solution time will not change significantly.

Procedure-A and Procedure-B are presented in detail in Section III.A and Section III.B respectively. The detailed SCED mathematical models are presented in Section IV.

### A. Procedure-A: SCED with RTCA

Fig. 2 illustrates the flowchart of the proposed Procedure-A that mimics existing industrial practice. As shown in Fig. 2, Procedure-A consists of the four steps listed below.
1) Monitor the system and determine the system status,
2) Perform base-case power flow and RTCA to determine the active network constraints for SCED,
3) Run SCED with the network constraints identified in 2),
4) Evaluate the SCED solution.

With the data collected from remote terminal units and local control centers, state estimation is performed to determine the system status in step 1. Then, base-case power flow and RTCA are conducted sequentially in step 2, which determines the active network constraints enforced in SCED for system secure operations. Then, SCED is solved and the optimal solution is reported to operators. As SCED uses the simplified DC power flow model, it is very important to ensure the AC feasibility of SCED solutions. The last step of Procedure-A evaluates the SCED solution by rerunning base-case power flow and contingency analysis with the updated generators' outputs; the contingency list used in this step is the same with step 2.

As EMS procedures are for real-time applications, the four steps of Procedure-A will repeat consecutively over time, which also holds for Procedure-B. In practice, SCED solutions must be approved by system operators before the unit dispatch signals can be sent out. It is possible that the AC power flow in the last step fails to converge, in which case the system operator would probably not approve the SCED solution. However, this rarely happens and is out of the scope of this paper.

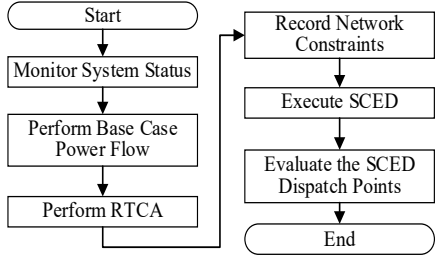

Fig. 2. Flowchart of the Proposed Procedure-A.

**Determination of Network Constraints**

Due to large-scale feature of real power systems, it is impractical to model all contingencies with all elements monitored in SCED. Thus, generation reserve is proposed and modeled in SCED aiming to have backup electric power for handling contingencies and enhancing reliability. However, reserve deliverability cannot be guaranteed due to congestion. Thus, SCED enforces a list of selected network constraints. With a limited number of network constraints, SCED can be solved within a short timeframe and be employed in real-time.

Network constraints can be divided into two categories, actual base-case network constraints and potential contingency-case network constraints. Each base-case network constraint contains three items: monitored branch $k$, initial flow $P_{k0}$ on branch $k$ in the base case, and the long-term normal limit $LimitA_k$. Each contingency-case network constraint contains four items: contingency $c$, monitored branch $k$, initial flow $P_{k,c,0}$ on branch $k$ under contingency $c$, and the short-term emergency limit $LimitC_{kc}$ under contingency $c$.

The actual thermal limits for transmission elements are for apparent power in the unit of MVA. However, SCED is performed to adjust generation only for active power in the unit of MW. Similarly, the branch thermal limits used in SCED are also for active power. Since the change in reactive power is insignificant in a short SCED period, branch MW limits used in SCED can be approximately derived by assuming the reactive power flows do not change in a short look-ahead period. Then, the $LimitA_k$ for base-case network constraint and the $LimitC_{kc}$ for contingency-case network constraint can be calculated by (1) and (2) respectively,

$$LimitA_k = \sqrt{RateA_k^2 - (\max(|Q_{k,f}|, |Q_{k,t}|))^2} \quad (1)$$

$$LimitC_{kc} = \sqrt{RateC_k^2 - (\max(|Q_{kc,f}|, |Q_{kc,t}|))^2} \quad (2)$$

where $Q_{k,f}$ and $Q_{k,t}$ denote the reactive powers on branch $k$ flowing out of from-bus and to-bus in base case respectively; $Q_{kc,f}$ and $Q_{kc,t}$ denote the reactive power on branch $k$ flowing out of from-bus and to-bus under contingency $c$ respectively.

When RTCA results are not available, branch emergency limit $LimitC_k$ can be approximately calculated by (3) which assumes reactive power does not change due to an outage.

$$LimitC_k = \sqrt{RateC_k^2 - (\max(|Q_{k,f}|, |Q_{k,t}|))^2} \quad (3)$$

To reduce computational complexity, only a small subset of branches will be monitored. Active network constraints are determined by comparing branch loading level with the tolerance $Pct$ for base case or the tolerance $PctC$ for contingency cases. Tolerances $Pct$ and $PctC$ are pre-defined percentages.

The branch loading level $LL_k$ for base case and $LL_{kc}$ for contingency case are defined in (4) and (5) respectively,

$$LL_k = \max(|S_{k,f}|, |S_{k,t}|)/RateA_k \quad (4)$$

$$LL_{kc} = \max(|S_{kc,f}|, |S_{kc,t}|)/RateC_k \quad (5)$$

where $S_{k,f}$ and $S_{k,t}$ denote the complex powers on branch $k$ flowing out of from-bus and to-bus in base case respectively; $S_{kc,f}$ and $S_{kc,t}$ denote the complex power on branch $k$ flowing out of from-bus and to-bus under contingency $c$ respectively.

Therefore, branch $k$ will be monitored in the base case if its loading level $LL_k$ is greater than $Pct$. Similarly, it will be monitored under contingency $c$ if the associated loading level $LL_{kc}$ is greater than $PctC$. The monitored branch constraints are referred to as active network constraints. An active network constraint is referred to as a critical network constraint if its initial flow exceeds the limit. A contingency is called an active contingency if it causes one or multiple active network constraints. Similarly, a critical contingency is a contingency that would cause at least one critical network constraint.

When $Pct$ is set to 1, only the congested branches and overloaded branches will be monitored in the base case. Similarly, when $PctC$ is set to 1, only the potentially congested branches and overloaded branches under contingency will be monitored for the associated contingencies. Typically, both $Pct$ and $PctC$ are set to unit as only congested lines and overloaded lines would be monitored. In practice, system operators can manually adjust the tolerances and lower the software reported branch limit by multiplying it with a percentage number such as 95% for system reliability concern at their own discretion.

The initial branch flows for base-case and contingency-case network constraints are determined by (6) and (7) respectively,

$$P_{k0} = sign(P_{k,f}) \cdot \max(|P_{k,f}|, |P_{k,t}|) \quad (6)$$

$$P_{kc0} = sign(P_{kc,f}) \cdot \max(|P_{kc,f}|, |P_{kc,t}|) \quad (7)$$

where $P_{k,f}$ and $P_{k,t}$ denote the active powers on branch $k$ flowing out of from-bus and to-bus in the base case respectively; $P_{kc,f}$ and $P_{kc,t}$ denote the active powers on branch $k$ flowing out of from-bus and to-bus under contingency $c$ respectively.

**Representation of Transmission Losses**

System losses typically accounts for 1% to 4% of the total demand. It would be impractical if all losses are compensated at a single slack bus. Thus, losses should be properly modeled in SCED. Losses can either be modeled as virtual loads or be calculated with loss coefficients. In this paper, virtual loads are used to represent transmission losses. As listed below, multiple methods are available to convert the losses in an AC model into virtual loads in a DC model:
- Assign losses to load buses,
- Assign losses to generator buses,
- Assign loss on each branch to the receiving bus,
- Assign loss on each branch to the sending bus,
- Assign loss on each branch evenly to its two buses.

In this paper, the loss on each branch is evenly distributed to the two buses that are connected to that branch and is modeled as fixed virtual loads. SCED runs consecutively every few minutes in real-time and typically adjusts the system generation slightly for online units only; the change in losses in a very short SCED period is not significant. Therefore, the losses are assumed to remain the same during a short SCED look-ahead period in this paper.

## B. Procedure-B: SCED with CTS-based RTCA

As stated in Section II, CTS can enhance system reliability by reducing post-contingency violation. Thus, in this section, Procedure-B is proposed to enhance Procedure-A by considering the benefits provided by CTS. Procedure-B can substantially relieve network congestion and significantly reduce congestion cost as compared to Procedure-A. To distinguish the regular SCED in Procedure-A, the SCED that implicitly models CTS in Procedure-B is referred to as enhanced SCED or E-SCED. The flowchart of the proposed Procedure-B is illustrated in Fig. 3 and it shows Procedure-B consists of six steps:

1) Monitor the system and determine the system status,
2) Perform base-case power flow and RTCA,
3) Perform CTS on critical contingencies only and identify beneficial switching actions for each critical contingency,
4) Update the thermal limits for critical network constraints,
5) Run E-SCED and obtain a new set of dispatch points,
6) Evaluate the E-SCED solution.

The first two steps of Procedure-B are the same with Procedure-A. The third step of Procedure-B is to perform CTS on critical contingencies only and identify the CTS actions that can reduce post-contingency violations, which aims to relieve network congestions with CTS. In step 4, the limits of critical network constraints are updated with (8). For a critical contingency $c$, if the switching solution can reduce the total violation while no individual violation is worse off, then, the pseudo limit of the associated constraint can be calculated by (8).

$$LimitC_{kc} = \sqrt{PRateC_{kc}^2 - (\max(|Q_{kc,f}|, |Q_{kc,t}|))^2} \quad (8)$$

where $PRateC_{kc}$ denotes the pseudo MVA limit of branch $k$ under contingency $c$ and is defined in the equation below,

$$PRateC_{kc} = RateC_k + v_{kc}P_{kc,CTS} \quad (9)$$

where $v_{kc}$ denotes the violation on branch $k$ under contingency $c$, and $P_{kc,CTS}$ denotes the percent violation reduction with CTS and is calculated by (10),

$$P_{kc,CTS} = (v_{kc} - v_{kc,CTS})/v_{kc} \quad (10)$$

where $v_{kc,CTS}$ is the flow violation on branch $k$ under contingency $c$ in the post-switching situation.

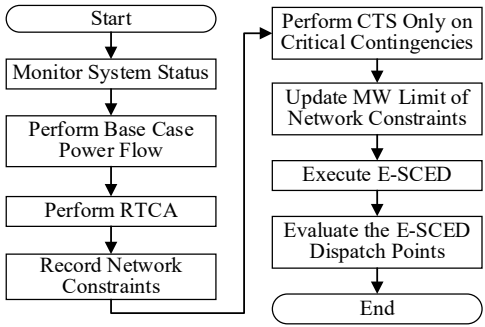

Fig. 3. Flowchart of the Proposed Procedure-B.

In step 5, E-SCED executes with the updated pseudo limits of critical network constraints and then obtains new generation points for dispatchable units. The last step of Procedure-B evaluates the E-SCED solutions. With the updated generations in the post-E-SCED stage, RTCA is performed on the same list of contingencies with step 2 and it may report violations as higher pseudo limits are used in E-SCED rather than the actual limits. Those violations would probably be on the same elements under the same contingencies as reported by RTCA in step 2. However, those violations will not be a concern; because the CTS solutions that are determined in step 3 in the pre-E-SCED stage can handle the same violations even if the system condition has changed in the post-E-SCED stage.

Procedure-A is a traditional EMS procedure that is consistent with industry practices. It contributes to this paper by disclosing the details about connection between full AC model based RTCA and simplified DC model based SCED. In addition, Procedure-A is implemented in this paper to (i) validate the effectiveness of the proposed SCED models in an AC framework, (ii) provide a basis to gauge the proposed Procedure-B with CTS in the EMS and measure the benefits by comparing the results of Procedure-A and Procedure-B.

## IV. SCED MODEL

Five SCED models are proposed in this section. The proposed SCED models co-optimize energy and reserve while enforcing physical restrictions and reliability requirements. Load shedding is included to handle infeasibility and prevent SCED from terminating without reporting any information.

### A. Unit Cost Curve

ISOs including PJM, Midcontinent ISO, and New York ISO require units to submit incremental offers that are represented by MW quantity and price pairs [34]-[36]. For instance, PJM accepts up to 10 price-quantity segments. There are two types of energy offers: block cost curve and slope cost curve.

Fig. 4 illustrates the generator block incremental cost curve. The lengths of the three segments are $P_{g,s1}$, $(P_{g,s2} - P_{g,s1})$, and $(P_{g,s3} - P_{g,s2})$ while the associated constant costs are $C_1$, $C_2$, and $C_3$ respectively. $p_{g1}$, $p_{g2}$, and $p_{g3}$ denote the net MW outputs on segments 1, 2, and 3 respectively. The optimal SCED solution will not schedule any power on a segment if any other segments with lower prices are not entirely selected since the objective is to minimize the total cost.

Fig. 5 illustrates the unit slope incremental cost curve. Obviously, the costs of the second and third segments are not constant, which creates non-linearity in the objective function. As non-linearity may create computational issues, linearization is desired. A slope segment can be divided into several sub-segments with the same length and thus, it can be approximately represented by a series of block sub-segments. Fig. 6 illustrates the linearization of a slope segment. The procedure for linearizing the second slope segment is presented below.

1) determine the number of sub-segments $nSS$ for a given slope segment. When minimum price increment $\Delta C$ is set, the nearest integer of $(C_2 - C_1)/\Delta C$ may be chosen to be $nSS$.
2) calculate the sub-segment breadth with (11),

$$\Delta s = (P_{g,s2} - P_{g,s1})/nSS \quad (11)$$

where $\Delta s$ denotes the actual breadth of each sub-segment.
3) the cost for each sub-segment can be determined by (12).

$$C_i = C_1 + (i - 0.5)\frac{(C_2 - C_1)}{(P_{g,s2} - P_{g,s1})}\Delta s, \ i = 1..nSS \quad (12)$$

With the above procedure, a slope segment can be converted into a series of small block segments. Therefore, a slope cost curve can be approximately represented by a block cost curve, which eliminates non-linearity in the model. Typically, for both block cost curves and slope cost curves, the first segment is flat, and it corresponds to the unit economic minimum generation $P_{g,s1}$ and the no-load cost $C_1$. Slope cost curves are converted into block cost curves before solving SCED.

### B. Objective Function

The proposed SCED models share the same objective function with existing industry practices, which is to minimize the total cost including operating energy cost and reserve cost. The objective function used in this paper is shown below,

$$\min EC + RC + C_{sh} \sum_{d \in D} p_{dsh} \quad (13)$$

where $EC$ denotes the energy cost and $RC$ denotes the reserve cost, and they are defined in (14) and (15) respectively.

$$EC = \sum_{g \in GD} \sum_{i=1}^{NS_g} p_{g,i} C_{g,i} \quad (14)$$
$$RC = \sum_{g \in G} sr_g CSR_g \quad (15)$$

### C. Constraints

The constraints of a basic power transfer distribution factor (PTDF)-based SCED model consist of (16)-(35). They can be divided into five categories: power balance constraints, load shedding constraints, generation constraints, reserve constraints, and network constraints. They are introduced below, as well as the alternative constraints.

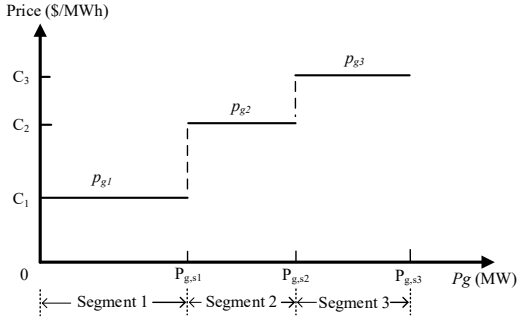
Fig. 4. Block cost curve of generator $g$.

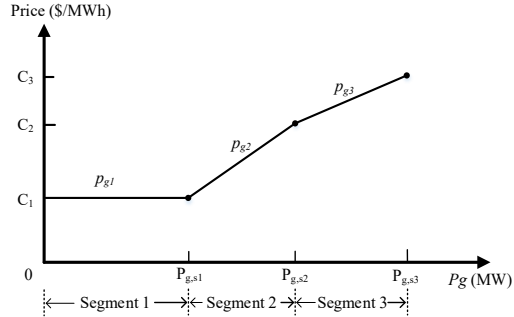
Fig. 5. Slope cost curve of generator $g$.

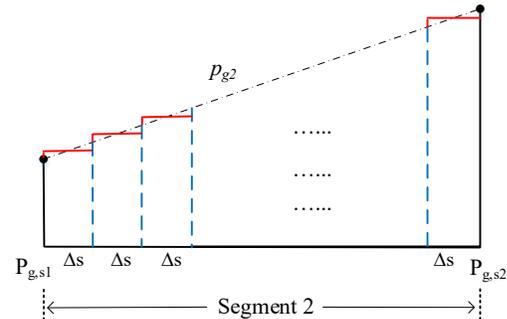
Fig. 6. Illustration of linearization of a slope segment.

**Power Balance Constraints**

For PTDF-based SCED, a single system-wide power balance constraint is sufficient to ensure the power balance, which is enforced in (16). However, a list of nodal power balance constraints must be used in $B$-$\theta$ based SCED models, which will be introduced later in this section.

$$\sum_{g \in G} p_g = \sum_{d \in D}(P_d - p_{d,sh}) \quad (16)$$

**Load Shedding Constraints**

In a practical power system, negative load may be used in the EMS to model small miscellaneous generations or tie line flows. It is not reasonable to shed loads that are used to represent generations or transfer flows with other neighboring systems. It is also not right to shed virtual loads that are used to represent losses. Thus, load shedding is not allowed for negative loads and virtual loads in this work, which is guaranteed by (17). The shedded load cannot exceed the actual positive demand, which is enforced by (18).

$$p_{d,sh} = 0 \;,\; d \in \{DN, DV\} \quad (17)$$
$$0 \leq p_{d,sh} \leq P_d \;,\; d \in \{D\} \backslash \{DN, DV\} \quad (18)$$

**Generation Constraints**

Some generators such as self-scheduling units may not be available for dispatch. Thus, those generators' outputs are fixed in SCED, which is expressed by (19). Dispatchable units typically offer incremental cost curves that consist of one or multiple pairs of price and quantity. Equation (20) ensures that a unit's generation equals the summation of the outputs on all segments. Constraint (21) guarantees that the power scheduled for each segment will not exceed the associated segment breadth. Constraint (22) enforces units' outputs to be within their upper limits and lower limits. Generators' energy ramping limit is modeled in (23).

$$p_g = P_{g0} \;,\; g \in (G - GD) \quad (19)$$
$$p_g = \sum_{i=1}^{NS_g} p_{g,i} \;,\; g \in GD \quad (20)$$
$$0 \leq p_{g,i} \leq BS_{g,i} \;,\; g \in GD \quad (21)$$
$$P_{g,min} \leq p_g \leq P_{g,max} \;,\; g \in G \quad (22)$$
$$-MRR_g T_{ED} \leq p_g - P_{g0} \leq MRR_g T_{ED} \;,\; g \in G \quad (23)$$

**Reserve Constraints**

The spinning reserve that an online unit can provide is subject to its spinning ramping capability, which is expressed in (24). Note that the ramping limits for energy re-dispatch and reserve deployment for the same unit may be different. Constraint (25) ensures that the sum of a unit's output and reserve does not exceed its maximum limit. In other words, reserve is also restricted by unit's available capacity in addition to ramping limit. As defined in (26), the "largest generator" rule is used for the reserve requirement, which ensures the reserve is sufficient to cover any of loss of a single generator.

$$0 \leq sr_g \leq SRR_g T_{SR} \;,\; \forall g \in G \quad (24)$$
$$p_g + sr_g \leq P_{g,max} \;,\; \forall g \in G \quad (25)$$
$$\sum_{m \in G} sr_m \geq p_g + sr_g \;,\; \forall g \in G \quad (26)$$

**Network Constraints**

Though reserve is modeled in SCED, network congestions may limit reserve deliverability. Thus, it is necessary to model active network constraints in SCED. Branch thermal limits for the base case and contingency cases are enforced in (27) and (28) respectively. The monitor sets can be different for base case and different contingency cases.

Due to concerns regarding voltage stability and transient stability, the total transfer power of the ties connecting two areas is not allowed to exceed a specific limit which is referred to as interface limit or transfer limit. In SCED, stability limit can be modeled by adding restriction on the sum of flows on branches that form an interface. The interface constraints for base case and contingency cases are represented by (29) and (30) respectively. The interface limit under contingency may be different with the limit in the base case especially when the contingency element is a branch of that interface.

Equations (31) and (32) calculate branch flows taking the effects of generation re-dispatch, load shedding, and load fluctuation into account. Note that (27) and (28) are only for branches in the monitor sets and (29) and (30) are only for critical interfaces; however, (31) and (32) are for branches in the monitor sets and branches forming the critical interfaces. Constraint (32) involves both LODF and outage transfer distribution factor (OTDF). The flow on contingency branch $c$ is forced to be zero in (33).

$$-LimitA_k \le p_k \le LimitA_k, k \in KM(0) \quad (27)$$
$$-LimitC_k \le p_{k,c} \le LimitC_k, k \in \text{KM}(c), c \in C \quad (28)$$
$$\sum_{k \in KI(i)} p_k \le Limit_i, i \in IM(0) \quad (29)$$
$$\sum_{k \in KI(i)} p_{k,c} \le Limit_{ic}, i \in IM(c), c \in C \quad (30)$$
$$p_k = P_{k0} + \sum_{n \in N}(PTDF_{n,k}(\Delta p_{gn} - \Delta p_{dn})), k \in \{KM(0), IKM(0)\} \quad (31)$$
$$p_{k,c} = P_{k0} + LODF_{k,c}P_{c0} + \sum_{n \in N} OTDF_{n,k,c}(\Delta p_{gn} - \Delta p_{dn}),$$
$$k \in \{KM(c), IKM(c)\} \backslash \{c\}, c \in C \quad (32)$$
$$p_{k,c} = 0, k \in \{c\}, c \in C \quad (33)$$

where $\Delta p_{gn}$ denotes change in generation at bus $n$ and $\Delta p_{dn}$ denotes change in load at bus $n$; and they are defined below,

$$\Delta p_{gn} = \sum_{g \in G(n)}(p_g - P_{g0}) \quad (34)$$
$$\Delta p_{dn} = \sum_{d \in D(n)}(P_d - P_{d0} - p_{d,sh}) \quad (35)$$

**Alternative Constraints**

Constraints (16)-(35) along with the objective function (13) form a basic SCED mathematical model. Some of those constraints can be replaced with alternative constraints. Enhancement and adjustment can be made to improve this model.

Adding $P_{g0}$ to each term in (23) would reformulate it to (36), which shares the same form with (22). By simply taking the minimum of upper limits as the new upper limit and using the maximum of lower limits as the new lower limit, constraints (36) and (22) can be combined as (37). In other words, constraints (22) and (23) can be replaced by one single constraint (37), which would reduce the number of constraints and increase the performance in terms of computational time.

$$P_{g0} - MRR_g T_{ED} \le p_g \le P_{g0} + MRR_g T_{ED}, g \in G \quad (36)$$
$$\max\{P_{g0} - MRR_g T_{ED}, P_{g,min}\} \le p_g \le \min\{P_{g0} + MRR_g T_{ED}, P_{g,max}\}, g \in G \quad (37)$$

When branch reactive power flow under contingency is available, the branch emergency limits can be customized for different contingencies. In other words, $LimitC_{kc}$ that is calculated by (2) should replace $LimitC_k$ in (28). Thus, (28) can be converted into the constraint shown below.

$$-LimitC_{k,c} \le p_{k,c} \le LimitC_{k,c}, k \in \text{KM}(c), c \in C \quad (38)$$

Similarly, if the initial branch flow $P_{k,c,0}$ under contingency $c$ is available from contingency analysis, the SCED model would be more accurate by replacing $P_{k0} + LODF_{k,c}P_{c0}$ with $P_{k,c,0}$ in (32). Then, (32) can be replaced by (39).

$$p_{k,c} = P_{k,c,0} + \sum_{n \in N}(OTDF_{n,k,c}(\Delta p_{gn} - \Delta p_{dn})), k \in \{KM(c), IKM(c)\} \backslash \{c\}, c \in C \quad (39)$$

If the initial branch flow for both the base case and contingency cases are not available, then, incremental PTDF based equations (31) and (32) can be replaced by cold-start PTDF based equations (40) and (41) respectively.

$$p_k = \sum_{n \in N}(PTDF_{n,k}(\sum_{g \in G(n)} p_g + \sum_{d \in D(n)}(p_{d,sh} - P_d))), k \in \{KM(0), IKM(0)\} \quad (40)$$
$$p_{k,c} = \sum_{n \in N}(OTDF_{n,k,c}(\sum_{g \in G(n)} p_g + \sum_{d \in D(n)}(p_{d,sh} - P_d))),$$
$$k \in \{KM(c), IKM(c)\} \backslash \{c\}, c \in C \quad (41)$$

Instead of using PTDF formulation, branch flow can be calculated using $B$-$\theta$ formulation, which are defined in (42) and (43). Power flow equations must be modeled for all branches with $B$-$\theta$ formulation because there are mutual effects between phase angles $\theta$ of all buses. It is worth mentioning that only the flows on branches of interests need to be calculated with PTDF formulation. With $B$-$\theta$ formulation based SCED, system-wide power balance constraint (16) should be replaced with nodal power balance constraints (44) and (45).

$$p_k = (\delta_{n(k-)} - \delta_{n(k+)} + \alpha_k)/x_k, k \in K \quad (42)$$
$$p_{k,c} = (\delta_{n(k-),c} - \delta_{n(k+),c} + \alpha_k)/x_k, k \in K \backslash \{c\}, c \in C \quad (43)$$
$$\sum_{g \in G(n)} p_g + \sum_{k \in K(n+)} p_k - \sum_{k \in K(n-)} p_k = \sum_{d \in D(n)}(P_d - p_{d,sh}), n \in N \quad (44)$$
$$\sum_{g \in G(n)} p_g + \sum_{k \in K(n+)} p_{k,c} - \sum_{k \in K(n-)} p_{k,c} = \sum_{d \in D(n)}(P_d - p_{d,sh}), n \in N, c \in C \quad (45)$$

*D. Models*

Based on different power flow formulations and the availability of network flow information, five different SCED models are proposed in this paper. They are listed below:
- SCED-M1: hot-start PTDF based SCED model,
- SCED-M2: warm-start PTDF based SCED model,
- SCED-M3: cold-start PTDF based SCED model,
- SCED-M4: hot-start $B$-$\theta$ based SCED model,
- SCED-M5: cold-start $B$-$\theta$ based SCED model.

A SCED model that takes the system initial condition information into account is called a hot-start model. Similarly, a SCED model with limited information or no information of the system condition is called a warm-start model or a cold-start model respectively.

The constraints enforced in different SCED models are listed in Table I. The first three SCED models use PTDF power flow formulation while the last two SCED models are based on $B$-$\theta$ power flow formulation. The five proposed SCED models share a common set of constraints. Moreover, the three PTDF based SCED models share the same system-wide power balance constraint while the two $B$-$\theta$ based SCED models share the same nodal power balance constraints.

In this paper, the proposed SCED models are preventive SCED that requires the system solution to be within limits under any single element outage without the need of generation re-dispatch in the post-contingency situation. To convert branch flow and thermal limit information from the full AC power flow model to the simplified DC power flow model, it is assumed that reactive power does not change in a short SCED period.

All five SCED models define the same variables including unit output $p_g$, unit output on each cost curve segment $p_{g,i}$, base case branch flow $p_k$, contingency case branch flow $p_{kc}$, unit spinning reserve $sr_g$, and shedded load $p_{d,sh}$. In addition, SCED-M4 and SCED-M5 also define base case phase angle $\delta_n$ and contingency case phase angle $\delta_{nc}$ due to the use of the $B$-$\theta$ power flow formulation. The most important decision variables for all five SCED models are the generator dispatch points $p_g$, which constitutes the SCED signals that would be sent from system control center to generator control centers.

Table I Constraints of multiple SCED models

| | Power balance constraints | Network constraints | Shared constraints |
|---|---|---|---|
| SCED-M1 | (16) | (31), (34)-(35), (38)-(39) | (17)-(21), (24)-(27), (29)-(30), (33), (37) |
| SCED-M2 | | (28), (31)-(32), (34)-(35) | |
| SCED-M3 | | (28), (40)-(41) | |
| SCED M4 | (44)-(45) | (38), (42)-(43) | |
| SCED M5 | | (28), (42)-(43) | |

The initial branch flow $P_{k,c,0}$ under contingency and the emergency limit $LimitC_{kc}$ are available to SCED-M1, while SCED-M2 uses $P_{k,0}$ and line outage distribution factor (LODF) to calculate branch flow under contingency and uses $LimitC_k$ as the emergency limits. A cold-start branch flow formulation is used in SCED-M3 rather than the incremental formulations used in SCED-M1 and SCED-M2.

Since SCED is built upon the linearized DC power flow model, apart from model approximation, the accuracy of flow calculation also depends on variation of net injection including generation. SCED-M3 uses total generator's output to calculate branch flow; however, incremental models only need to consider the change in generation when calculating branch flow. Thus, the model precision of SCED-M3 would be less than SCED-M1 and SCED-M2.

SCED-M4 and SCED-M5 are based on traditional $B$-$\theta$ power flow model. The difference between them is that SCED-M4 uses $LimitC_{kc}$ as the customized emergency limits for different contingency cases while SCED-M5 uses the same emergency limit $LimitC_k$ for all contingency cases.

The five proposed SCED models differ in the formulation of network constraints that are based on different forms of power flow equations and network flow information. Their AC performances will be tested and compared using Procedure-A in the Part-II of this paper and the SCED model with the best performance will be selected for demonstrating the effectiveness of Procedure-B. In practical situations, not all required information for the best SCED model will always be available; thus, another important information obtained from the model comparisons is the accuracy level of each SCED model, which can help operators properly evaluate the effects of different SCED models and address the issues associated with model imperfection.

*E. Market Implication*

In addition to providing updated dispatch solutions in real-time, SCED is also used to determine energy market solutions including locational marginal prices (LMPs). Thus, it is very important to analyze the effect of integrating CTS into SCED on market results. LMP, load payment, generator revenue, generator cost, generator rent, congestion cost, and congestion revenue are presented and analyzed in this paper.

Locational marginal pricing is a market mechanism that is used to clear wholesale energy markets. The nodal LMP reflects the least cost of supplying the next increment load at that specific bus while meeting all physical and reliability constraints. LMP consists of three components: energy component, congestion component, and losses component. If a system has a network with infinite capacity and no losses, LMP would be the same over the entire system. However, in reality, losses cannot be avoided and network congestion typically exists. In this paper, losses are represented by virtual loads in the proposed SCED models. Thus, the losses component of LMP is ignored in this paper and the nodal LMP can then be calculated by the following equation,

$$LMP_n = LMP_s + LMP_{cg,n}, \ n \in N \quad (46)$$

where congestion component $LMP_{cg,n}$ is calculated below,

$$LMP_{cg,n} = \sum_{k \in KM(0)} PTDF_{n,k}(F_k^+ - F_k^-) +$$
$$\sum_{c \in C} \sum_{k \in KM(c)} OTDF_{n,k,c}(F_{k,c}^+ - F_{k,c}^-) +$$
$$\sum_{i \in IM(0)} \sum_{k \in KI(i)} PTDF_{n,k}(F_i^+ - F_i^-) +$$
$$\sum_{c \in C} \sum_{i \in IM(0)} \sum_{k \in KI(i)} OTDF_{n,k,c}(F_{i,c}^+ - F_{i,c}^-), n \in N \quad (47)$$

where $F_k^+$ and $F_k^-$ denote dual variables of thermal limit constraints of branch $k$ in the base case; $F_{kc}^+$ and $F_{kc}^-$ are dual variables of thermal limit constraints of branch $k$ under contingency $c$; $F_i^+$ and $F_i^-$ denote dual variables of stability limit constraints of interface $i$ in the base case; $F_{kc}^+$ and $F_{kc}^-$ are dual variables of stability limit constraints of interface $i$ under contingency $c$.

*Proof for* (46)

Duality theory is used to prove (46). Since this proof is only for deriving nodal LMP expression (46), the constraints and variables that are not of interest are ignored for simplicity, as well as the objective function. The constraints and associated dual variables of interest are listed below.

$$\sum_{n \in N} PI_n = 0 \quad (LMP_s) \quad (48)$$
$$\sum_{g \in G(n)} p_g + \sum_{d \in D(n)} p_{d,shed} - PI_n = \sum_{d \in D(n)} P_d, n \in N \quad (LMP_n) \quad (49)$$
$$P_{k0} + \sum_{n \in N}(PTDF_{n,k}(PI_n - PI_{n,0})) \leq LimitA_k, k \in KM(0) \quad (F_k^+) \quad (50)$$
$$-P_{k0} - \sum_{n \in N}(PTDF_{n,k}(PI_n - PI_{n,0})) \leq LimitA_k, k \in KM(0) \quad (F_k^-) \quad (51)$$
$$P_{k,c,0} + \sum_{n \in N}(OTDF_{n,k}(PI_n - PI_{n,0})) \leq LimitC_{kc}, k \in KM(c), c \in C \quad (F_{k,c}^+) \quad (52)$$
$$-P_{k,c,0} - \sum_{n \in N}(OTDF_{n,k}(PI_n - PI_{n,0})) \leq LimitC_{kc}, k \in KM(c), c \in C \quad (F_{k,c}^-) \quad (53)$$
$$\sum_{k \in KI(i)}(P_{k0} + \sum_{n \in N}(PTDF_{n,k}(PI_n - PI_{n,0}))) \leq Limit_i, i \in IM(0) \quad (F_i^+) \quad (54)$$
$$-\sum_{k \in KI(i)}(P_{k0} + \sum_{n \in N}(PTDF_{n,k}(PI_n - PI_{n,0}))) \leq Limit_i, i \in IM(0) \quad (F_i^-) \quad (55)$$
$$\sum_{k \in KI(i)}(P_{k,c,0} + \sum_{n \in N}(OTDF_{n,k}(PI_n - PI_{n,0}))) \leq Limit_{ic}, i \in IM(c)), c \in C \quad (F_{i,c}^+) \quad (56)$$
$$-\sum_{k \in KI(i)}(P_{k,c,0} + \sum_{n \in N}(OTDF_{n,k}(PI_n - PI_{n,0}))) \leq Limit_{ic}, i \in IM(c)), c \in C \quad (F_{i,c}^-) \quad (57)$$

Constraint (48) ensures system-wide power balance while (49) is for nodal power balance. Constraint (48) is redundant for $B$-$\theta$ based SCED while (49) is not needed for PTDF based SCED; however, they are listed here just for deriving the relationship between nodal LMP and system LMP. Constraints (50)-(53) show that the system is subject to branch thermal limit. Moreover, power systems are also restricted by interface limits, which is guaranteed by (54)-(57). Variables $LMP_s$, $LMP_n$, $F_k^+$, $F_k^-$ $F_{k,c}^+$, $F_{k,c}^-$, $F_i^+$, $F_i^-$, $F_{i,c}^+$, and $F_{i,c}^-$ are dual variables of constraints (48)-(57) in this primal problem. Note that variable $PI_n$ denotes net power injection at bus $n$, which can be either positive or non-positive; then, it is unconstrained in

IEEE Transactions on Power Systems                                                                                                                      10the primal problem. Therefore, the associated constraint in the dual problem is an equality constraint as expressed in (58).

$$LMP_n - LMP_s + \sum_{k \in KM(0)} PTDF_{n,k}(F_k^+ - F_k^-) + \sum_{c \in C} \sum_{k \in KM(c)} OTDF_{n,k,c}(F_{k,c}^+ - F_{k,c}^-) + \sum_{i \in IM(0)} \sum_{k \in KI(i)} PTDF_{n,k}(F_i^+ - F_i^-) + \sum_{c \in C} \sum_{i \in IM(0)} \sum_{k \in KI(i)} OTDF_{n,k,c}(F_{i,c}^+ - F_{i,c}^-) = 0, \ n \in N \quad (58)$$

where $LMP_s$ and $LMP_n$ are unconstrained; $F_k^+$, $F_k^-$, $F_{k,c}^+$, $F_{k,c}^-$, $F_i^+$, $F_i^-$, $F_{i,c}^+$, and $F_{i,c}^-$ are non-positive. Then, (46) can be easily derived by reformatting (58).

As defined in (59), average LMP over the entire system is proposed to analyze the effect of modeling CTS in SCED on LMP. Similarly, average congestion LMP is defined in (60).

$$AvgLMP = \sum_{n \in N} LMP_n /nb, \ n \in N \quad (59)$$
$$AvgLMP_{cg} = \sum_{n \in N} LMP_{cg,n}/nb, \ n \in N \quad (60)$$

Load payment is calculated by (61) and generator revenue is determined by (62). Equation (63) calculates the generator cost which is part of the objective function. Generator rent is defined in (64). Note that in this paper, generator rent only accounts for energy and does not include reserve rent. Congestion revenue, which is used to fund the financial transmission rights markets, is the difference between load payment and generator revenue, and it can be calculated by (65).

$$LdPaymt = \sum_{n \in N}(LMP_n(\sum_{d \in D(n)} P_d)) \quad (61)$$
$$GenRvn = \sum_{n \in N}(LMP_n(\sum_{g \in G(n)} p_g)) \quad (62)$$
$$GenCost = \sum_{g \in GD} \sum_{i=1}^{NS_g} p_{g,i} C_{g,i} \quad (63)$$
$$GenRent = GenRvn - GenCost \quad (64)$$
$$CngstRvn = LdPaymt - GenRvn \quad (65)$$

In addition to congestion revenue, congestion cost is also proposed to measure the degree of network congestion. Congestion cost is defined as (66),

$$CngstCost = TCost1 - TCost2 \quad (66)$$

where $TCost1$ denotes the total cost of either an E-SCED or a SCED and $TCost2$ denotes the total cost for the same E-SCED or SCED problem but without any network constraints.

Thus, the congestion cost reduction achieved by E-SCED as compared to a traditional SCED can be calculated by (67).

$$CCR_{CTS} = CngstCost_{ESCED} - CngstCost_{SCED} \quad (67)$$

## V. CONCLUSIONS

The authors' prior work has demonstrated that CTS can significantly reduce post-contingency violations and illustrated CTS can enhance system reliability with flexible transmission. However, the economic benefits that may be achieved with CTS are not studied. Thus, this paper investigated how much economic benefits can be translated from the reliability benefits provided by CTS.

Two EMS procedures are proposed in this two-part paper. Procedure-A, which connects AC based RTCA and DC based SCED, is consistent with industrial practice. Extended upon Procedure-A, Procedure-B is proposed to enable E-SCED to take advantage of the flexibility in transmission networks in a practical way by using the proposed pseudo limit of network constraints rather than the actual limit for a traditional SCED.

Part-I includes a comprehensive literature review on RTCA, RT SCED, and CTS, and presents the proposed methodology. The case studies conducted in Part-II demonstrate that 1) the proposed Procedure-A can efficiently eliminate the actual violations and the potential post-contingency overloads identified by RTCA, 2) with consideration of CTS in E-SCED, the proposed Procedure-B can achieve significant economic benefits by relieving the potential post-contingency network congestion, and 3) the CTS solution that reduces branch overloads for a contingency in the pre-SCED stage can also provide benefits for the same contingency in the post-SCED stage.

## REFERENCES

[1] K. W. Hedman, S. S. Oren, and R. P. O'Neill, "A review of transmission switching and network topology optimization," *IEEE PES General Meeting*, Detroit, MI, Jul. 2011.
[2] California ISO, "Minimum effective threshold report," Mar. 2010, [Online]. Available: http://www.caiso.com/274c/274ce77df630.pdf.
[3] PJM, "Manual 3: transmission operations - section 5: index and operating procedures for PJM RTO operation," Jun. 2010, [Online]. Available: http://www.pjm.com/~/media/training/nerc-certifications/m03v37-transmission-operations.ashx.
[4] ISO-NE, "ISO New England operating procedure No. 19: transmission operations," [Online]. Available: http://www.iso-ne.com/rules_proceds/operating/isone/op19/op19_rto_final.pdf.
[5] X. Li, P. Balasubramanian, M. Abdi-Khorsand, A. S. Korad, and K. W. Hedman, "Effect of topology control on system reliability: TVA test case," *CIGRE Grid of the Future Symposium*, Houston, TX, Oct. 2014.
[6] X. Li, P. Balasubramanian, M. Sahraei-Ardakani, M. Abdi-Khorsand, K. W. Hedman, R. Podmore, "Real-Time Contingency Analysis with Correct Transmission Switching," *IEEE Trans. Power Syst.*, vol. 32, no. 4, pp. 2604-2617, Jul. 2017.
[7] M. Sahraei-Ardakani, X. Li, P. Balasubramanian, K. W. Hedman, and M. Abdi-Khorsand, "Real-time contingency analysis with transmission switching on real power system data," *IEEE Trans. Power Syst.*, vol. 31, no. 3, pp. 2501-2502, May 2016.
[8] X. Li, P. Balasubramanian, K. W. Hedman, "A Data-driven Heuristic for Corrective Transmission Switching," *North American Power Symposium (NAPS)*, Denver, CO, Sep. 2016.
[9] P. Balasubramanian, M. Sahraei-Ardakani, X. Li, K. W Hedman, "Towards Smart Corrective Switching: Analysis and Advancement of PJM's Switching Solutions," *IET Generation, Transmission, and Distribution*, vol. 10, no. 8, pp. 1984-1992, May 2016.
[10] W. Shao and V. Vittal. "Corrective switching algorithm for relieving overloads and voltage violations," *IEEE Trans. Power Syst.*, vol. 20, no. 4, pp. 1877-1885, Nov. 2005.
[11] E. B. Fisher, R. P. O'Neill, and M. C. Ferris, "Optimal Transmission Switching," *IEEE Trans. Power Syst.*, vol. 23, no. 3, pp. 1346-1355, Aug. 2008.
[12] K. W. Hedman, R. P. O'Neill, E. B. Fisher, and S. S. Oren, "Optimal transmission switching with contingency analysis," *IEEE Trans. Power Syst.*, vol. 24, no. 3, pp. 1577–1586, Aug. 2009.
[13] S. Fattahi, J. Lavaei, and A. Atamtürk, "A bound strengthening method for optimal transmission switching in power systems," *IEEE Trans. Power Syst.*, vol. 34, no. 1, pp. 280-291, Jan. 2019.
[14] P. A. Ruiz, J. M. Foster, A. Rudkevich, and M. C. Caramanis, "Tractable Transmission Topology Control Using Sensitivity Analysis," *IEEE Trans. Power Syst.*, vol. 27, no. 3, pp. 1550-1559, 2012.
[15] J. D. Fuller, R. Ramasra, and A. Cha, "Fast Heuristics for Transmission-line Switching," *IEEE Trans. Power Syst.*, vol. 27, no. 3 pp. 1377-1386, Aug. 2012.
[16] J. Han and A. Papavasiliou, "Congestion management through topological corrections: A case study of Central Western Europe," *Energy Policy*, vol. 86, pp. 470-482, Nov. 2015.
[17] G. Granelli, M. Montagna, F. Zanellini, P. Bresesti, R. Vailati, and M. Innorta, "Optimal Network Reconfiguration for Congestion Management by Deterministic and Genetic Algorithms," *Electric Power Systems Research*, vol. 76, no. 6-7, pp. 549-556, Apr. 2006.
[18] M. Abdi-Khorsand and K. W. Hedman, "Day-ahead corrective transmission topology control," *IEEE PES General Meeting*, Washington D.C., USA, Jul. 2014.
[19] A. Khodaei and M. Shahidehpour, "Transmission Switching in Security-constrained Unit Commitment," *IEEE Trans. Power Syst.*, vol. 25, no. 4, pp. 1937-1945, Nov. 2010.
[20] K. W. Hedman, M. C. Ferris, R. P. O'Neill, E. B Fisher, and S. S. Oren, "Co-optimization of generation unit commitment and transmission switching with N-1 reliability," *IEEE Trans. Power Syst.,* vol. 25, no. 2, pp. 1052-1063, May 2010.
the primal problem. Therefore, the associated constraint in the dual problem is an equality constraint as expressed in (58).

$$LMP_n - LMP_s + \sum_{k \in KM(0)} PTDF_{n,k}(F_k^+ - F_k^-) + \sum_{c \in C} \sum_{k \in KM(c)} OTDF_{n,k,c}(F_{k,c}^+ - F_{k,c}^-) + \sum_{i \in IM(0)} \sum_{k \in KI(i)} PTDF_{n,k}(F_i^+ - F_i^-) + \sum_{c \in C} \sum_{i \in IM(0)} \sum_{k \in KI(i)} OTDF_{n,k,c}(F_{i,c}^+ - F_{i,c}^-) = 0, \ n \in N \quad (58)$$

where $LMP_s$ and $LMP_n$ are unconstrained; $F_k^+$, $F_k^-$, $F_{k,c}^+$, $F_{k,c}^-$, $F_i^+$, $F_i^-$, $F_{i,c}^+$, and $F_{i,c}^-$ are non-positive. Then, (46) can be easily derived by reformatting (58).

As defined in (59), average LMP over the entire system is proposed to analyze the effect of modeling CTS in SCED on LMP. Similarly, average congestion LMP is defined in (60).

$$AvgLMP = \sum_{n \in N} LMP_n /nb, \ n \in N \quad (59)$$

$$AvgLMP_{cg} = \sum_{n \in N} LMP_{cg,n}/nb, \ n \in N \quad (60)$$

Load payment is calculated by (61) and generator revenue is determined by (62). Equation (63) calculates the generator cost which is part of the objective function. Generator rent is defined in (64). Note that in this paper, generator rent only accounts for energy and does not include reserve rent. Congestion revenue, which is used to fund the financial transmission rights markets, is the difference between load payment and generator revenue, and it can be calculated by (65).

$$LdPaymt = \sum_{n \in N}(LMP_n(\sum_{d \in D(n)} P_d)) \quad (61)$$

$$GenRvn = \sum_{n \in N}(LMP_n(\sum_{g \in G(n)} p_g)) \quad (62)$$

$$GenCost = \sum_{g \in GD} \sum_{i=1}^{NS_g} p_{g,i} C_{g,i} \quad (63)$$

$$GenRent = GenRvn - GenCost \quad (64)$$

$$CngstRvn = LdPaymt - GenRvn \quad (65)$$

In addition to congestion revenue, congestion cost is also proposed to measure the degree of network congestion. Congestion cost is defined as (66),

$$CngstCost = TCost1 - TCost2 \quad (66)$$

where $TCost1$ denotes the total cost of either an E-SCED or a SCED and $TCost2$ denotes the total cost for the same E-SCED or SCED problem but without any network constraints.

Thus, the congestion cost reduction achieved by E-SCED as compared to a traditional SCED can be calculated by (67).

$$CCR_{CTS} = CngstCost_{ESCED} - CngstCost_{SCED} \quad (67)$$

## V. CONCLUSIONS

The authors' prior work has demonstrated that CTS can significantly reduce post-contingency violations and illustrated CTS can enhance system reliability with flexible transmission. However, the economic benefits that may be achieved with CTS are not studied. Thus, this paper investigated how much economic benefits can be translated from the reliability benefits provided by CTS.

Two EMS procedures are proposed in this two-part paper. Procedure-A, which connects AC based RTCA and DC based SCED, is consistent with industrial practice. Extended upon Procedure-A, Procedure-B is proposed to enable E-SCED to take advantage of the flexibility in transmission networks in a practical way by using the proposed pseudo limit of network constraints rather than the actual limit for a traditional SCED.

Part-I includes a comprehensive literature review on RTCA, RT SCED, and CTS, and presents the proposed methodology. The case studies conducted in Part-II demonstrate that 1) the proposed Procedure-A can efficiently eliminate the actual violations and the potential post-contingency overloads identified by RTCA, 2) with consideration of CTS in E-SCED, the proposed Procedure-B can achieve significant economic benefits by relieving the potential post-contingency network congestion, and 3) the CTS solution that reduces branch overloads for a contingency in the pre-SCED stage can also provide benefits for the same contingency in the post-SCED stage.

**Xingpeng Li** (S'13–M'18) received the B.S. degree in electrical engineering from Shandong University, Jinan, China, in 2010, and the M.S. degree in electrical engineering from Zhejiang University, Hangzhou, China, in 2013. He received the M.S. degree in industrial engineering and the Ph.D. degree in electrical engineering from Arizona State University, Tempe, AZ, USA, in 2016 and 2017 respectively.

Currently, he is an Assistant Professor in the Department of Electrical and Computer Engineering at the University of Houston. He previously worked for ISO New England, Holyoke, MA, USA, and PJM Interconnection, Audubon, PA, USA. Before joining the University of Houston, he was a Senior Application Engineer for ABB, San Jose, CA, USA. His research interests include power system operations and optimization, energy management system, energy markets, microgrids, and grid integration of renewable energy sources.

**Kory W. Hedman** (S'05–M'10–SM'16) received a B.S. degree in electrical engineering and a B.S. degree in economics from the University of Washington, Seattle, WA, USA, in 2004, an M.S. degree in economics and an M.S. degree in electrical engineering from Iowa State University, Ames, IA, USA, in 2006 and 2007, respectively, and the M.S. and Ph.D. degrees in operations research from the University of California—Berkeley, Berkeley, CA, USA, in 2008 and 2010, respectively.

He is currently an Associate Professor with the School of Electrical, Computer, and Energy Engineering, Arizona State University, Tempe, AZ, USA. His research interests include power systems operations and planning, market design, power system economics, renewable energy, and operations research. He was the recipient of the Presidential Early Career Award for Scientists and Engineers from the U.S. President Barack H. Obama in 2017.